# A SIMPLE PROOF OF KAIJSER'S UNIQUE ERGODICITY RESULT FOR HIDDEN MARKOV $\alpha$-CHAINS

By Fred Kochman and Jim Reeds

*Center for Communications Research*

According to a 1975 result of T. Kaijser, if some nonvanishing product of hidden Markov model (HMM) stepping matrices is subrectangular, and the underlying chain is aperiodic, the corresponding $\alpha$-chain has a unique invariant limiting measure $\lambda$.

Here the $\alpha$-chain $\{\alpha_n\} = \{(\alpha_{ni})\}$ is given by

$$\alpha_{ni} = P(X_n = i | Y_n, Y_{n-1}, \ldots),$$

where $\{(X_n, Y_n)\}$ is a finite state HMM with unobserved Markov chain component $\{X_n\}$ and observed output component $\{Y_n\}$. This defines $\{\alpha_n\}$ as a stochastic process taking values in the probability simplex. It is not hard to see that $\{\alpha_n\}$ is itself a Markov chain. The stepping matrices $M(y) = (M(y)_{ij})$ give the probability that $(X_n, Y_n) = (j, y)$, conditional on $X_{n-1} = i$. A matrix is said to be subrectangular if the locations of its nonzero entries forms a cartesian product of a set of row indices and a set of column indices.

Kaijser's result is based on an application of the Furstenberg–Kesten theory to the random matrix products $M(Y_1)M(Y_2)\cdots M(Y_n)$. In this paper we prove a slightly stronger form of Kaijser's theorem with a simpler argument, exploiting the theory of e chains.

**1. Introduction.** In 1975 Kaijser [9] gave a simple sufficient condition for the uniqueness of the invariant measure for the so-called $\alpha$-chain, a certain weak Feller chain with compact state space arising in the study of an arbitrary finite state hidden Markov model. This provided an elegant partial answer to a question posed by David Blackwell in 1957 [3]. (We follow Blackwell in using $\alpha_n$ to denote the state of the $\alpha$-chain at time $n$. Kaijser calls it $Z_n$.) The transition behavior of a finite state hidden Markov model (HMM) and of its associated $\alpha$-chain can be specified by a finite collection of substochastic matrices, the *stepping matrices*; and probability calculations with HMMs involve, at least conceptually, lengthy matrix products of stepping matrices. Accordingly, Kaijser's analysis utilized the Furstenburg–Kesten theory [6] of random matrix products. However, by exploiting the theory of e-chains, in particular Theorem 18.4.4 of [11], we are able to give a simpler proof of a result slightly stronger than that in [9].

Briefly, an HMM [2] consists of a pair of stochastic processes, $\{X_n\}$ and $\{Y_n\}$, taking values in finite sets $\mathcal{X}$ and $\mathcal{Y}$, such that $\{X_n\}$ is a Markov chain and each $Y_n$ is a probabilistic function of $(X_{n-1}, X_n)$. In modeling applications [5, 10], the "observable" marginal process $\{Y_n\}$ is "output" from the "hidden" process $\{X_n\}$.

The transition structure of $(X_n, Y_n)$ can be specified by the *stepping matrices* $M(y) = (M(y)_{ij})$, given by

$$M(y)_{ij} = P((X_{n+1}, Y_{n+1}) = (j, y)|X_n = i);$$

their sum $M = \sum_y M(y)$ is equal to the transition matrix of the Markov chain $\{X_n\}$.

Let $\Delta$ be the finite-dimensional simplex of probability measures on $\mathcal{X}$, and provisionally set $\alpha_0 \in \Delta$ to be the marginal distribution of $X_0$ and, for $n \geq 0$, set $\alpha_n \in \Delta$ to be the conditional distribution of $X_n$ given $\{Y_t : 1 \leq t \leq n\}$. It can be shown [3, 9], that $\{\alpha_n\}$ is a Markov chain with $\Delta$ as its (continuous) state space.

Blackwell [3] first studied the $\alpha$-chain and posed the question of when its transition law has a unique invariant measure. His partial answer is based on contractivity hypotheses far stronger than Kaijser's hypotheses, or ours. For Blackwell, $\alpha_n$ is the conditional distribution of $X_n$ given the infinite past $\{Y_t : -\infty < t \leq n\}$, where now $\{X_n\}$ is assumed to be stationary, so $\{\alpha_n\}$ is stationary as well. This $\{\alpha_n\}$ is again a Markov chain, with the same transition law as the provisional $\{\alpha_n\}$ defined above. (Blackwell's motivation for studying the $\alpha$-chain is a formula expressing the entropy of $\{Y_n\}$ in terms of the distribution of his version of $\alpha_n$.) By starting the $\alpha$-chain in the infinite past, so it is in effect born stationary, Blackwell avoids questions of convergence to a limiting distribution. But by allowing the $X$ chain to start at finite time $t = 0$, with arbitrary distribution, Kaijser opens the possibility that the $\alpha$-chain could be nonstationary, which, in turn, raises the additional question about whether the finite-time distributions of $\alpha_n$ converge to a stationary limit measure. These two versions of the $\alpha$-chain have the same transition mechanism, but usually different initial or marginal distributions on $\Delta$.

For our result, we allow $\{\alpha_n\}$ to have any initial distribution on $\Delta$ at time $t = 0$, but with the same transition law as above. Of course, this destroys the original motivating interpretation of $\alpha_n$ as a conditional distribution of $X_n$; but it is the (unchanged) transition mechanism whose properties are of primary interest to us, rather than any particular realization of the chain.

We now prepare to state Kaijser's result. Following Kaijser, call a nonnegative matrix $(D_{ij})$ *subrectangular* if the set of subscript pairs $(i,j)$ with $D_{ij} > 0$ forms a Cartesian product, that is, if there exist sets $R$ and $C$ of row and column subscripts so that $D_{ij} > 0$ if and only if $(i,j) \in R \times C$. Call the matrix $M$ (and the chain $\{X_n\}$) *irreducible* if for every pair of subscripts $(i,j)$, there is some $k$ [whose value may depend on $(i,j)$] for which $(M^k)_{ij} > 0$. If there is a single value of $k$ such that for all $(i,j)$ we have $(M^k)_{ij} > 0$, the matrix $M$ (and the chain $\{X_n\}$) are said to be *aperiodic*.

With this terminology, the result is (slightly paraphrased):

THEOREM A ([9]). *Suppose the transition matrix $M$ is aperiodic. Suppose some nonzero product of stepping matrices is subrectangular. Then the probability distribution of $\alpha_n$ converges weakly to a unique limit measure, independent of the initial distribution for $\alpha_0$.*

Kaijser's argument is along the following lines. A path for the chain $\{\alpha_n\}$, starting from $\alpha_0$, can be written as

$$\alpha_0, \frac{\alpha_0 M(Y_1)}{\alpha_0 M(Y_1) e}, \frac{\alpha_0 M(Y_1) M(Y_2)}{\alpha_0 M(Y_1) M(Y_2) e}, \ldots, \frac{\alpha_0 M(Y_1) M(Y_2) \cdots M(Y_k)}{\alpha_0 M(Y_1) M(Y_2) \cdots M(Y_k) e}, \ldots,$$

where $\{Y_n\}$ is the marginal process defined above and $e$ is the column vector of all 1's. Since the sequence $\{M(Y_n)\}$ of matrices is itself a stochastic process, Kaijser was able to cleverly adapt methods of the Furstenberg–Kesten theory to the present subject.

However, we have a different line of argument, based on the theory of e-chains, which we think is ultimately easier to understand. We replace the subrectangularity condition with the following, which we will show to be weaker. Let $\mathcal{M}$ be the set of stepping matrices, let $\mathcal{M}^*$ be the set of all finite products of elements of $\mathcal{M}$, and let $\mathcal{C} = \mathbb{R}^+ \mathcal{M}^*$ be the cone on $\mathcal{M}^*$, that is, all positive scalar multiples of elements of $\mathcal{M}^*$. Our condition is that the closure, $\overline{\mathcal{C}}$, should contain a matrix of rank 1.

A very brief sketch of our argument is as follows.

First, in Theorem 1 we show our key technical result, that for an arbitrary transition matrix $M$ and arbitrary decomposition into stepping matrices, $\{\alpha_n\}$ is an *e-chain*, in the sense of [11], page 144. Then we exploit the associated limit theory. Namely, given our rank 1 hypothesis, when the matrix $M$ is irreducible, we show that the state space $\Delta$ possesses a *topologically reachable point $v$*, in the sense of [11], page 455. Further, if $M$ is also aperiodic, then $v$ must be *topologically aperiodic* in the sense of [11], page 459, as well. Since $\Delta$ is compact, Theorem 18.4.4 of [11], page 460, immediately applies, yielding the proof of our main result:

THEOREM 2. *Let the matrix $M$ be irreducible. Suppose there exists a rank 1 element of $\overline{\mathcal{C}}$. Then $T_{\mathcal{M}}$ has a unique stationary distribution $\lambda$, and $(T_{\mathcal{M}}^*)^n \mu \longrightarrow \lambda$ weakly in Cesaro mean, for each probability measure $\mu$ on $\Delta$. If, in addition, $M$ is aperiodic, then also $(T_{\mathcal{M}}^*)^n \mu \longrightarrow \lambda$ weakly, for each probability measure $\mu$ on $\Delta$.*

[Here $T_{\mathcal{M}}$ denotes the Markov operator on $C(\Delta)$ associated with the $\alpha$-chain.]

Finally, we derive Kaijser's Theorem A from our Theorem 2, by showing that if his subrectanglarity hypothesis is true, so is our rank 1 hypothesis.

We conclude the paper with three calculations. The first shows that aperiodicity of $M$, by itself, does not imply uniqueness of the stationary measure for the $\alpha$-chain. Another shows that in Kaijser's theorem the condition of aperiodicity cannot be replaced with that of irreducibility. The third shows that an example of Kaijser's, while not satisfying the conditions of his result (Theorem A), does satisfy the conditions of ours (Theorem 2).

We now address the relation of this work to "random systems with complete connections" [4] and Chapter 2 of [8] and "place dependent random iterated function systems" (IFS) [1]. The $\alpha$-chain seems—ignoring technicalities—to fall under the scope of these theories, so one might suppose Kaijser's theorem followed as a corollary of standard IFS results. We have, however, been unable to derive Kaijser's results this way. Our main obstacle is that the state update functions for $\alpha$-chains are only defined, in general, on open dense subsets of $\Delta$ and need not extend continuously to all of $\Delta$; nor do they seem to satisfy the conventional contractivity or mean contractivity hypotheses imposed in the IFS literature. Since we ultimately rely on the classical Perron theorem for aperiodic matrices, we too are exploiting a kind of contractivity; the difference seems to be that it enters at a later stage of the argument.

In common with Kaijser's argument, ours does exploit the special role played by matrix products in $\alpha$-chain calculations. We find it striking how smoothly the theory of e-chains may be applied without clutter to the matrix product formulation, once the necessary ground work is completed. There does not seem to be any easy analogue of this matrix product structure in the generic IFS example.

**2. Notation, formulae.** Let $\mathcal{X}$ and $\mathcal{Y}$ be finite sets, with $s = |\mathcal{X}|$, and let $\Delta \subset \mathbb{R}^s$ be the simplex of probability measures on $\mathcal{X}$. Let $e$ be the $s$-long column vector of all 1's. Let $C(\Delta)$ be the space of all continuous real-valued functions on $\Delta$, with the sup-norm topology. Let $\mathcal{P}(\Delta)$ denote the probability measures on $\Delta$, equipped with with the weak topology. Let $\mathcal{M} = \{M(y) : y \in \mathcal{Y}\}$ be a family of nonnegative matrices whose sum, $M = \sum_y M(y)$, is a Markov transition matrix; so the *stepping matrices* $M(y)$ are substochastic. Let $\mathcal{M}^*$ be the set of finite products of elements of $\mathcal{M}$. As a convenience, we will use notations like $\mathbf{y} = y_1, y_2, \ldots, y_n \in \mathcal{Y}^n$ to denote finite sequences of elements of $\mathcal{Y}$. If $\mathbf{y} \in \mathcal{M}^n$, let $|\mathbf{y}| = n$ denote its length. If $\mathbf{y} = y_1, y_2, \ldots, y_n \in \mathcal{Y}^n$, we will let $M(\mathbf{y})$ be shorthand for the matrix product $M(y_1)M(y_2)\cdots M(y_n)$. We use the term *word* to refer (indiscriminately) to tuples $\mathbf{y}$ or matrix products $M(\mathbf{y})$.

For each $\mu \in \mathcal{P}(\Delta)$, the $\alpha$-chain may be concisely defined as follows: Pick $\alpha_0$ according to $\mu$; conditional on $\alpha_0$, pick $\{Y_n\}$ so that $P(Y_1 = y_1, \ldots, Y_n = y_n) = \alpha_0 M(y_0) \cdots M(y_n)e$, and then conditional on $\alpha_0$ and $Y_1, \ldots, Y_n$, define $\{\alpha_n\}$ to satisfy the conditionally certain recursion

$$\alpha_n = \frac{\alpha_{n-1} M(Y_n)}{\alpha_{n-1} M(Y_n)e},$$

which is to say,

$$\alpha_n = \frac{\alpha_0 M(Y_1) \cdots M(Y_n)}{\alpha_0 M(Y_1) \cdots M(Y_n)e}.$$

Though $\{Y_n\}$ is not generally Markov of any finite order, it is known [3, 9] that $\{\alpha_n\}$ is a Markov chain on $\Delta$ whose transition law is given by the transition kernel

$$P(\alpha, A) = {\sum_y}' \alpha M(y)e,$$

where the sum extends over all $y \in \mathcal{Y}$ such that

$$\frac{\alpha M(y)}{\alpha M(y)e} \in A.$$

The chain is weak Feller: its Markov operator

$$T_{\mathcal{M}} : C(\Delta) \to C(\Delta)$$

is given by the formula

$$(T_{\mathcal{M}} f)(\alpha) = \sum_y (\alpha M(y)e) f\left(\frac{\alpha M(y)}{\alpha M(y)e}\right),$$

where any term with $\alpha M(y)e = 0$ is set equal to 0. For later use, we record the telescoped $n$-step transition formulas

(1) $$P^n(\alpha, A) = {\sum_\mathbf{y}}' \alpha M(\mathbf{y})e,$$

where the sum extends over all $\mathbf{y} \in \mathcal{M}^n$ such that $\alpha M(\mathbf{y})/\alpha M(\mathbf{y})e \in A$, and

(2) $$(T_{\mathcal{M}}^n f)(\alpha) = \sum_{\mathbf{y} \in \mathcal{M}^n} (\alpha M(\mathbf{y})e) f\left(\frac{\alpha M(\mathbf{y})}{\alpha M(\mathbf{y})e}\right).$$

Since $\Delta$ is compact, it is immediate that at least one $T_{\mathcal{M}}$-invariant probability law exists; part of what is at issue is whether there is more than one. If there is a unique $T_{\mathcal{M}}$-invariant probability law, we say the $\alpha$-chain is *uniquely ergodic*.

**3. The $\alpha$-chain is an e-chain.** A weak Feller chain with compact state space $S$ is an *e-chain* if its operator $T$ on $C(S)$ is such that, for each $f \in C(S)$, the set of functions $\{T^n f : n > 0\}$ is equicontinuous. We show that every $\alpha$-chain is an e-chain:

THEOREM 1. *Let $M$ be any transition matrix with any stepping decomposition $\mathcal{M}$. Then the corresponding $\alpha$-chain is an e-chain.*

PROOF. By the Arzela–Ascoli theorem, since $\Delta$ is compact, it suffices to show that, for given $f \in C(\Delta)$, the set $\{T_{\mathcal{M}}^n f\}$ is relatively compact. To this end, we will construct a compact set $K_f \subset C(\Delta)$ for which $\{T_{\mathcal{M}}^n f\} \subseteq K_f$.

Let $h \in \Delta$ be some fixed probability distribution on $\mathcal{X}$ for which $h_i > 0$ for all $i \in \mathcal{X}$. Let $\mathcal{W} = \{V = (v_{ij}) : v_{ij} \geq 0, hVe = 1\}$ be the $s \times s$ matrices with nonnegative entries obeying the linear constraint $hVe = 1$.

Given $f$, we define a continuous function $g(\alpha, V) = \alpha V e f(\alpha V / \alpha V e)$. The function is defined in the first instance for those $(\alpha, V) \in \Delta \times \mathcal{W}$ for which $\alpha V e \neq 0$, and because $f$ is bounded, $g$ has a unique continuous extension to all of $\Delta \times \mathcal{W}$: take $g(\alpha, V) = 0$ when $\alpha V e = 0$. Let $K_f \subset C(\Delta)$ be the set of all functions of $\alpha$ obtained by integrating $g(\alpha, V)$ with respect to all probability measures on the compact set $\mathcal{W}$, that is, all functions of form

$$\alpha \mapsto Eg(\alpha, V)$$

for random elements $V$ in $\mathcal{W}$. Thus, $K_f$ is the closed convex hull of the compact set of all the functions of $\alpha$ obtained from $g(\alpha, V)$ by holding $V$ fixed. Hence, $K_f$ is also compact.

For given $n$, pick a random element $\mathbf{w} \in \mathcal{Y}^n$ with probability $hM(\mathbf{w})e$ and set

$$V_n = \frac{M(\mathbf{w})}{hM(\mathbf{w})e},$$

which, with probability 1, is a matrix in $\mathcal{W}$. Then, referring to (2), we see that

$$(T_{\mathcal{M}}^n f)(\alpha) = Eg(\alpha, V_n),$$

so $T_{\mathcal{M}}^n f \in K_f$. Thus, $\{T_{\mathcal{M}}^n f\}$ is contained in $K_f$ and, hence, is equicontinuous. □

**4. Main result.** We now embark on the proof of our main result. Recall that $\mathcal{C}$ is the set of all positive scalar multiples of matrices in $\mathcal{M}^*$ and that its closure is $\overline{\mathcal{C}}$.

THEOREM 2. *Let $M$ and $\mathcal{M}$ be given. Suppose $M$ is irreducible and suppose that $\overline{\mathcal{C}}$ contains an element of rank $1$. Then $T_{\mathcal{M}}$ has a unique stationary distribution $\lambda$, and $(T_{\mathcal{M}}^*)^n \mu \longrightarrow \lambda$ weakly in Cesaro mean, for each $\mu \in \mathcal{P}(\Delta)$. Suppose, in addition, that $M$ is aperiodic. Then $(T_{\mathcal{M}}^*)^n \mu \longrightarrow \lambda$ weakly, for each $\mu \in \mathcal{P}(\Delta)$.*

PROOF. We will use the rank 1 element of $\overline{\mathcal{C}}$ to construct a topologically reachable point $v \in \Delta$. If $M$ is aperiodic, $v$ will also be topologically aperiodic. According to Theorem 18.4.4 of [11], page 460, the existence of such $v$, the fact that we are working with an e-chain, and the compactness of $\Delta$ together imply the stated results.

Suppose $R \in \overline{\mathcal{C}}$ has rank 1, so $R = uv$, where $u \neq 0$ is a nonnegative column vector and $v$ a nonnegative row vector which we may assume scaled so $ve = 1$. In particular, $v \in \Delta$. We will show that if $M$ is irreducible, then $v$ is topologically reachable, that is, for each $\alpha \in \Delta$ and each open set $\mathcal{O}$ containing $v$, there exists a $k > 0$ such that $P^k(\alpha, \mathcal{O}) > 0$.

For each $\alpha \in \Delta$, there is some word $M(\mathbf{z})$ such that $\alpha M(\mathbf{z}) u > 0$, as follows. There are certainly $i$ and $j$ with $\alpha_i > 0, u_j > 0$, and since for some $k$, we have $\sum_{\mathbf{z} \in \mathcal{M}^k} M(\mathbf{z})_{ij} = (M^k)_{ij} > 0$, we must have $\alpha M(\mathbf{z}) u > 0$ for some $\mathbf{z} \in \mathcal{M}^k$.

As a consequence, $\alpha M(\mathbf{z}) R$ is a nonzero multiple of $v$. But $R$ is a limit of rescaled words:

$$R = \lim_{n \to \infty} M(\mathbf{y}_n)/s_n$$

for some sequence of words $\mathbf{y}_n$ and reals $s_n > 0$. For all $n$ sufficiently large, $\alpha M(\mathbf{z}) M(\mathbf{y}_n) e > 0$, so

$$v = \frac{\alpha M(\mathbf{z}) R}{\alpha M(\mathbf{z}) R e} = \lim_{n \to \infty} \frac{\alpha M(\mathbf{z}) M(\mathbf{y}_n)}{\alpha M(\mathbf{z}) M(\mathbf{y}_n) e}.$$

This implies that, for any neighborhood $\mathcal{O}$ of $v$, for $n$ large enough, we must also have

$$\frac{\alpha M(\mathbf{z}) M(\mathbf{y}_n)}{\alpha M(\mathbf{z}) M(\mathbf{y}_n) e} \in \mathcal{O}.$$

Hence, referring to (1), $P^k(\alpha, \mathcal{O}) \geq \alpha M(\mathbf{z}) M(\mathbf{y}_n) e > 0$ for $k = |\mathbf{z}| + |\mathbf{y}_n|$, and so $v$ is topologically reachable.

If $M$ is also aperiodic, then a strengthening of this argument yields a positive lower bound on $P^k(\alpha, \mathcal{O})$ which is uniform in large $k$, showing that $v$ is a topologically aperiodic state. Let $\mathbf{y}_n$, $u$, $v$ and $R = uv$ be as above, and let $\pi$ be the stationary probability vector for $M$. Pick $\mathbf{z}$ so $\pi M(\mathbf{z})u > 0$, let $w = M(\mathbf{z})u/\pi M(\mathbf{z})u$, and define $S_n = M(\mathbf{z})M(\mathbf{y}_n)$. Then $\lim_{n\to\infty} S_n/\pi S_n e = wv$.

Let $\|\cdot\|$ denote the $l_1$ norm for row vectors and the induced operator norm for matrices acting on row vectors on the right, so for row vector $\alpha$ and matrix $T$ we have $|\alpha Te| \leq \|\alpha T\| \leq \|\alpha\|\|T\|$. Let $\mathcal{B}$ be the closed unit $l_1$ ball in $\mathbb{R}^s$. Then there exist matrices $T_n$ and scalars $\delta_n \geq 0$ so that

$$\frac{S_n}{\pi S_n e} = wv + \delta_n T_n,$$

with $\|T_n\| \leq 1$ and $\lim_{n\to\infty} \delta_n = 0$. Now pick $n$ so large that $\pi S_n e > 0$ and that $\delta_n$ is sufficiently small that both $\delta_n < 1/4$ and, for all $\beta \in \mathcal{B}$,

$$\frac{v + \sqrt{\delta_n}\beta}{1 + \sqrt{\delta_n}\beta e} \in \mathcal{O}.$$

Finally, let $t = 2\sqrt{\delta_n}\pi S_n e$ and let $m = |\mathbf{z}| + |\mathbf{y}_n|$.

Given all these choices, we claim that, for all $\alpha \in \Delta$,

(3) $$P^{k+m}(\alpha, \mathcal{O}) \geq \alpha M^k S_n e - t$$

for all $k \geq 0$. If so, since $M$ is aperiodic, $\alpha M^k \to \pi$ as $k \to \infty$, so

$$\liminf_{k\to\infty} P^{k+m}(\alpha, \mathcal{O}) \geq \pi S_n e - t$$

$$= (1 - 2\sqrt{\delta_n})\pi S_n e > 0.$$

Letting $\alpha = v$ shows, in particular, that $v$ is a topologically aperiodic state.

To prove (3), first assume $k = 0$. Since (3) is then trivially true if $\alpha S_n e \leq t$, we may assume $\alpha S_n e > t$. But in that case $S_n$ steps $\alpha$ into $\mathcal{O}$ as follows. Since

$$\alpha w + \delta_n \alpha T_n e = \frac{\alpha S_n e}{\pi S_n e} > 2\sqrt{\delta_n},$$

we get a lower bound on the scalar $\alpha w$:

$$\alpha w > 2\sqrt{\delta_n} - \delta_n \alpha T_n e > 2\sqrt{\delta_n} - \delta_n > \sqrt{\delta_n}.$$

Let $\beta = \sqrt{\delta_n}\,\alpha T_n/\alpha w$, so $\beta \in \mathcal{B}$. Then

$$\frac{\alpha S_n}{\alpha S_n e} = \frac{\alpha S_n}{\pi S_n e} \bigg/ \frac{\alpha S_n e}{\pi S_n e}$$

$$= \frac{\alpha wv + \delta_n \alpha T_n}{\alpha w + \delta_n \alpha T_n e}$$

$$= \frac{v + \sqrt{\delta_n}\beta}{1 + \sqrt{\delta_n}\beta e} \in \mathcal{O}.$$

Hence, the word $S_n \in \mathcal{M}^m$ steps $\alpha$ into $\mathcal{O}$. This implies $P^m(\alpha, \mathcal{O}) \geq \alpha S_n e$, verifying (3) when $k = 0$.

For $k > 0$, we have

$$P^{k+m}(\alpha, \mathcal{O}) = \sum_{|\mathbf{w}|=k} \alpha M(\mathbf{w}) e \, P^m\left(\frac{\alpha M(\mathbf{w})}{\alpha M(\mathbf{w})e}, \mathcal{O}\right)$$

$$\geq \sum_{|\mathbf{w}|=k} \alpha M(\mathbf{w}) e \left(\frac{\alpha M(\mathbf{w}) S_n e}{\alpha M(\mathbf{w})e} - t\right)$$

$$= \sum_{|\mathbf{w}|=k} \alpha M(\mathbf{w}) S_n e - \sum_{|\mathbf{w}|=k} \alpha M(\mathbf{w}) e t$$

$$= \alpha M^k S_n e - t,$$

concluding the verification of (3).

By Theorem 1, $\{\alpha_n\}$ is an e-chain; it is obviously *bounded in probability on average* in the sense of [11], page 285, since $\Delta$ is compact.

Hence, Theorem 18.4.4 of [11], page 460, applies, and our theorem follows.
$\square$

**5. Kaijser's result.** We are now in a position to derive Kaijser's theorem from our Theorem 2.

THEOREM A ([9]). *Suppose $M$ is aperiodic. If there is a nonvanishing subrectangular $M(\mathbf{y}) \in \mathcal{M}^*$, then there exists a unique $T_\mathcal{M}$-invariant probability measure $\lambda$, and for all $\mu \in \mathcal{P}(\Delta)$, we have $(T_\mathcal{M}^*)^n \mu \longrightarrow \lambda$ (weakly) as $n \to \infty$.*

PROOF. First, we find a nonvanishing subrectangular word $G = M(\mathbf{z})$ with a positive entry in its $(1,1)$ position. If $M(\mathbf{y})$ does not already have this property, we pick $(i,j)$ such that $M(\mathbf{y})_{ij} > 0$, and then, as in the proof of Theorem 2, find words $M(\mathbf{u})$ and $M(\mathbf{v})$ such that $M(\mathbf{u})_{1i} > 0$ and $M(\mathbf{v})_{j1} > 0$. Let $\mathbf{z} = \mathbf{uyv}$. The product of a subrectangular matrix and a nonnegative matrix is subrectangular, so $M(\mathbf{z}) = M(\mathbf{u})M(\mathbf{y})M(\mathbf{v})$ has the desired property.

Let $R$ and $C$ be the sets of row and column indices which specify where $G_{ij} > 0$. That is, $G_{ij} > 0$ if and only if $i \in R$ and $j \in C$. For notational convenience, pretend that $R = S_\mathrm{I} \cup S_\mathrm{II}$ and $C = S_\mathrm{I} \cup S_\mathrm{III}$, where $S_\mathrm{I}, S_\mathrm{II}, S_\mathrm{III}$ and $S_\mathrm{IV}$ are a partition of $\mathcal{X}$ into blocks of consecutive integers. That is to say, $G$ has block structure

$$G = \begin{pmatrix} A & 0 & B & 0 \\ C & 0 & D & 0 \\ 0 & 0 & 0 & 0 \\ 0 & 0 & 0 & 0 \end{pmatrix},$$

where all of the entries in blocks $A, B, C$ and $D$ are strictly positive, and the blocks on the diagonal are square. In particular, since $G_{11} > 0$, the upper left block $A$ is $k \times k$, for some $k \geq 1$.

Check by induction that, for $n \geq 2$,

$$G^n = \begin{pmatrix} A \\ C \\ 0 \\ 0 \end{pmatrix} A^{n-2} \begin{pmatrix} A & 0 & B & 0 \end{pmatrix}.$$

By the Perron theorem [7], page 502, a suitably rescaled version of $A^n$ has a limit:

$$\lim_{n \to \infty} \theta^n A^n = \bar{A},$$

where $\theta > 0$ is the reciprocal of the spectral radius of $A$, all elements of $\bar{A}$ are strictly positive, and $\bar{A}$ has rank 1.

Hence, for some sequence of scaling constants $s_n$,

$$\lim_{n \to \infty} G^n / s_n = \begin{pmatrix} A \\ C \\ 0 \\ 0 \end{pmatrix} \bar{A} \begin{pmatrix} A & 0 & B & 0 \end{pmatrix}$$

also has rank 1. Since each $G^n \in \mathcal{M}^*$, we have exhibited a rank 1 element of $\overline{\mathcal{C}}$ and the result then follows from Theorem 2. $\square$

**6. Three examples and a question.** First, we give an example showing that the assumption of aperiodicity, by itself, is not enough to guarantee unique ergodicity of the $\alpha$-chain. The matrices

$$M = \begin{pmatrix} 1/2 & 1/2 \\ 1/2 & 1/2 \end{pmatrix}, \qquad M(0) = \begin{pmatrix} 1/2 & 0 \\ 0 & 1/2 \end{pmatrix}, \qquad M(1) = \begin{pmatrix} 0 & 1/2 \\ 1/2 & 0 \end{pmatrix}$$

give rise to an $\alpha$-chain with the following simple description: $\alpha_n = (u, v) \in \Delta \subset \mathbb{R}^2$ moves with probability $1/2$ to $\alpha_{n+1} = (u, v)$ and with probability $1/2$ to $(v, u)$. Thus, $|u - v|$ is a nontrivial invariant for the $\alpha$-chain, which therefore has multiple stationary distributions. Examples of such include the following: the uniform distribution on $\Delta$, the point mass at $(1/2, 1/2)$ and, for each $0 < u < 1/2$, the measures assigning probability $1/2$ to each of $(u, 1-u)$ and $(1-u, u)$. (A similar example appears in [9].)

Next, we give an example showing that the assumption of aperiodicity cannot be replaced by irreducibility in Kaijser's result. The matrices

$$M = \begin{pmatrix} 0 & 1 \\ 1 & 0 \end{pmatrix}, \qquad M(0) = \begin{pmatrix} 0 & 0 \\ 1 & 0 \end{pmatrix}, \qquad M(1) = \begin{pmatrix} 0 & 1 \\ 0 & 0 \end{pmatrix}$$

specify an HMM satisfying the subrectangularity condition. $M$ is clearly irreducible but not aperiodic. If the starting measure $\mu$ puts mass 1 at $\alpha_0 = (x, 1-x)$, where $x \notin \{0, 1/2, 1\}$, the subsequence $\alpha_{2n}$ has one limit distribution [which puts mass $x$ at $(1,0)$ and mass $1-x$ at $(0,1)$] and the subsequence $\alpha_{2n+1}$ has a different limit distribution [which puts mass $x$ at $(0,1)$ and mass $1-x$ at $(1,0)$].

Finally, at the end of his paper Kaijser conjectures that if $p \neq q$, the $\alpha$-chain for the HMM with two stepping matrices

$$M(0) = \begin{pmatrix} p & 0 & 0 & 0 \\ 0 & 1/2 & 0 & 0 \\ 1/2 & 0 & 0 & 0 \\ 0 & 1/2 & 0 & 0 \end{pmatrix} \quad \text{and} \quad M(1) = \begin{pmatrix} 0 & 0 & q & 0 \\ 0 & 0 & 0 & 1/2 \\ 0 & 0 & 0 & 1/2 \\ 0 & 0 & 1/2 & 0 \end{pmatrix}$$

has a unique stationary distribution, even though there are no nonzero subrectangular words. This conjecture is true, as we now show.

Applying the method used in our proof of Theorem A, consider the rescaled limits of $M(0)^n$. It is easy to check that

$$M(0)^n = \begin{pmatrix} p^n & 0 & 0 & 0 \\ 0 & 1/2^n & 0 & 0 \\ p^n/2p & 0 & 0 & 0 \\ 0 & 1/2^n & 0 & 0 \end{pmatrix}.$$

So if $p > 1/2$,

$$\lim_{n \to \infty} \frac{M(0)^n}{e' M(0)^n e} = \frac{2p}{2p+1} \begin{pmatrix} 1 & 0 & 0 & 0 \\ 0 & 0 & 0 & 0 \\ 1/2p & 0 & 0 & 0 \\ 0 & 0 & 0 & 0 \end{pmatrix}$$

and if $p < 1/2$,

$$\lim_{n \to \infty} \frac{M(0)^n}{e' M(0)^n e} = \begin{pmatrix} 0 & 0 & 0 & 0 \\ 0 & 1/2 & 0 & 0 \\ 0 & 0 & 0 & 0 \\ 0 & 1/2 & 0 & 0 \end{pmatrix}.$$

In either case the limit has rank 1, so, by Theorem 2, the $\alpha$-chain has a unique invariant measure.

Thus, Kaijser's subrectangularity condition is sufficient but not necessary.

In light of our proof, as well as this example, one may ask the following question: Is the condition that $\overline{\mathcal{C}}$ contains a rank 1 matrix a necessary and sufficient condition for the $\alpha$-chain to have a unique invariant measure, when $M$ is irreducible?

**Acknowledgment.** The authors are grateful to a referee for suggestions which greatly improved the presentation.

CENTER FOR COMMUNICATIONS RESEARCH
PRINCETON, NEW JERSEY 08540
USA
E-MAIL: kochman@idaccr.org
       reeds@idaccr.org